\documentclass[13pt,reqno]{article}

\usepackage{mathptmx}
\usepackage{amssymb}
\usepackage{amsfonts}
\usepackage{amsthm}
\usepackage{amsmath}
\usepackage{graphicx}
\usepackage{shadow}
\usepackage[all]{xy}
\usepackage{enumerate}
\usepackage{makeidx}
 \def\R{\hbox{{\rm I}\kern-0.2em{\rm R}\kern0.2em}}%mathematical R for reals

 \def\a{\alpha}  
 
 \def\be{\begin{equation}} \def\ee{\end{equation}} 
  \def\p{\partial} \def\({\left(} \def\){\right)} 
 \def\[{\left[}
 \def\]{\right]}
 \def\bc{\begin{center}}
 \def\ec{\end{center}}

 \def\bc{\begin{center}}
\def\ec{\end{center}}

\begin{document}

\bc \Large{\bf Generalization of the double reduction theory }
 \ec

\medskip

\bc Ashfaque H. Bokhari, Ahmad Y. Dweik, F. D. Zaman\\

{Department of Mathematical Sciences, King Fahd University of
Petroleum and Minerals, Dhahran 31261, Saudi Arabia}\\
A. H. Kara\\
{School of Mathematics, University of the Witwatersrand, Wits 2050, South Africa.}\\
F. M. Mahomed\\

{School of Computation and Applied Mathematics, Center for
Differential Equations, Continuum Mechanics and Applications,
University of the Witwatersrand, Wits 2050,  South Africa}\\
\ec
%%%%%%%%%%%%%%%%%%%%%%%%%%%%%%%%%%%%%%%%%%%%%%%%%%%%%%%%%%%%%%%%%%%%%%%%%%%%%%%%%%%%%%%%%%%%
\begin{abstract}
In a recent work \cite{art:2007,art:2008} Sj\"{o}berg remarked
that generalization of the double reduction theory to partial
differential equations of higher dimensions is still an open
problem. In this note we have attempted to provide this
generalization to find invariant solution for a non linear system
of $q$th  order partial differential equations with $n$
independent and $m$ dependent variables provided that the non
linear system of partial differential equations admits a
nontrivial conserved form which has at least one associated
symmetry in every reduction. In order to give an application of the procedure we apply it to the  nonlinear (2 + 1) wave
equation for arbitrary function $f(u)$ and $g(u)$.

\end{abstract}
%%%%%%%%%%%%%%%%%%%%%%%%%%%%%%%%%%%%%%%%%%%%%%%%%%%%%%%%%%%%%%%%%%%%%%%%%%%%%%%%%%%%%%%%%%%%%
\section{Introduction}
Applying a Lie point or Lie-B\"acklund  symmetry generator to a
conserved vector provide either (1) Conservation law associated
with that symmetry or (2) Conservation law that may be trivial,
known already or new. A pioneering work in this direction was
published by Kara et. al \cite{art:2000,art:2002}. Sj\"{o}berg later showed that \cite{art:2007,art:2008} when the
generated conserved vector is null, i.e. the symmetry is
associated with the conserved vector (association defined as in
\cite{art:2000}), a double reduction is possible for PDEs with two
independent variables. In this double reduction the PDE of order
$q$ is reduced to an ODE of order $(q - 1)$. Thus the use of one
symmetry associated with a conservation law leads to two
reductions, \emph{the first being a reduction of the number of
independent variables} and \emph{the second being a reduction of
the order of the DE.}
Sj\"{o}berg also constructed the reduction formula for PDEs with
two independent variables which transform the conserved form of
the PDE to a reduced conserved form via an associated symmetry.
Application of this method to the linear heat, the BBM and the
sine-Gordon equation and a system of differential equations from
one dimensional gas dynamics are given \cite{art:2007}.
The double reduction theory says that  a PDE of order
 $q$ with two independent and $m$ dependent variables, which
admits a nontrivial conserved form that has at least one
associated symmetry can be reduced to an ODE of order  $(q-1)$.\\
 In her papers \cite{art:2007,art:2008} Sj\"{o}berg opines that
 generalizing the double reduction theory to PDEs of higher dimensions is still
an open problem and it is not clear how to overcome the problem
when not all derivatives of non-local variables are known
explicitly. Further calculations for higher dimensions are quite
tedious and cumbersome. There do not exist enough examples of
potential symmetries and symmetries with associated conservation
laws for higher dimensional PDEs so that the complexity of this
problem can be demonstrated. And much work is needed to generalize
(if possible) the theory to PDEs with more than two independent
variables.\\
In this article we discuss \emph{a generalization of the double
reduction theory }with $n$ independent variables by showing that a
non linear system of $q$th order PDEs with $n$ independent and $m$
dependent variables, which admits a nontrivial conserved form that
has at least one associated symmetry in every reduction from the
$n$ reductions(the first step of double reduction) can be reduced
to a non linear system of $(q-1)$th order ODEs.\\
In order to solve this we use two main steps:
(a) Generalize the reduction formula of  Sj\"{o}berg in
\cite{art:2007} from two independent variable to $n$ independent
variables and (b) prove that the conserved form of PDEs with $n$ independent
variables can be transformed to a reduced conserved form via an
associated symmetry. Finally we apply the generalized double reduction to
the  nonlinear (2 + 1) wave equation for arbitrary function $f(u)$
and $g(u)$ to obtain invariant solution.

%%%%%%%%%%%%%%%%%%%%%%%%%%%%%%%%%%%%%%%%%%%%%%%%%%%%%%%%%%%%%%%%%%%%%%%%%%%%%%%%%%%%%%%%%%%%%%%%%%%
\section {The Fundamental Theorem of double reduction}
\setcounter{equation}{0} Consider the $q$th-order system of
partial differential equations (PDEs) of $n$ independent variables
 $x = (x^{1}, x^{2}, ... , x^{n})$ and $m$ dependent variables $u = (u^{1}, u^{2},...,
 u^{m})$
 \begin{equation}\label{a1}
 \begin{array}{ll}
 E^{\a}(x,u,u_{(1)},...,u_{(q)})=0, & \a=1,...,m~,
\end{array}
\end{equation}
where $u_{(1)},u_{(2)},...,u_{(q)}$ denote the collections of all
first, second,..., $q$ th-order partial derivatives, i.e.,
$u_{i}^{\a}=D_{i}(u^{\a}),u_{ij}^{\a}=D_{j}D_{i}(u^{\a})$,...respectively,
with the total differentiation operator with respect to $x^{i}$
given by
\begin{equation}\label{a2}
\begin{array}{ll}
D_{i}=\frac{\p}{\p x^{i}}+u_{i}^{\a}\frac{\p}{\p
u^\a}+u_{ij}^\a\frac{\p}{\p u_{j}^\a}+...,&i=1,...,n~,
\end{array}
\end{equation}
in which the summation convention is used.\\
The following definitions are well-known (see, e.g. \cite{art:1982,art:1997,art:2000}).\\
\emph{The Lie-B\"acklund operator} is
\begin{equation}\label{a4}
\begin{array}{ll}
X=\xi^{i}\frac{\p}{\p {x^i}}+\eta^\a\frac{\p}{\p
u^\a}&\xi^i,\eta^\a\in \emph{A}~,
\end{array}
\end{equation}
where $A$ is the space of  \emph{differential functions}. The
operator (\ref{a4}) is an abbreviated form of the infinite formal
sum
\begin{equation}\label{a5}
\begin{array}{ll}
 X=\xi^{i}\frac{\p}{\p x^{i}}+\eta^{\a}\frac{\p}{\p u^{\a}}+\sum\limits_{s\geq 1 }\zeta_{i_{1}i_{2}...i_{s}}^{\a}\frac{\p}{\p u_{i_1 i_2...i_s}^\a},\\
\end{array}
\end{equation}
where the additional coefficients are determined uniquely by the
prolongation formulae,
\begin{equation}\label{a6}
\begin{array}{ll}
\zeta_{i}^\a=D_{i}(W^\a)+\xi^{j}u_{ij}^\a&\\
\zeta_{i_1...i_s}^\a=D_{i_1}...D_{i_s}(W^\a)+\xi^j
u_{ji_1...i_s}^\a,& s>1,
\end{array}
\end{equation}
in which $W^\a$ is the\emph{ Lie characteristic function,}
\begin{equation}\label{a7}
W^\a=\eta^\a-\xi^j u_{j}^{\a}.
\end{equation}
The $n$-tuple vector $T = (T^1,T^2,...,T^n)$,   $T^j \in A,$    $j
= 1,...,n $ is a \emph{\emph{conserved vector}} of (\ref{a1}) if
$T^i$ satisfies
\begin{equation}\label{a12}
D_iT^i\mid_{(2.1)}=0.
\end{equation}
A Lie-B\"{a}cklund symmetry generator $X$ of the form (\ref {a5}) is associated with a
conserved vector $T$ of the system (\ref{a1}) if $X$ and $T$
satisfy the relations
\begin{equation}\label{a19}
\begin{array}{cc}
$[$T^i,X$]$=X(T^i)+T^{i}D_k(\xi^k)-T^{k}D_{k}(\xi^i)=0,&i=1,...,n.
\end{array}
\end{equation}
\textbf{Theorem 2.1} \cite{art:1982,art:1997} Suppose that $X$ is
any Lie-B\"{a}cklund symmetry of (\ref{a1}) and \\$T^i,i=1,...,n$
are the components of conserved
vector of (\ref{a1}).Then
\begin{equation}\label{a20}
\begin{array}{cc}
{T^*}^i=$[$T^i,X$]$=X(T^i)+T^{i}D_j\xi^j-T^{j}D_{j}\xi^i,  &i=1,...,n.\\
\end{array}
\end{equation}
constitute the components of a conserved vector of (\ref{a1}) ,i.e.\\
$D_i {T^*}^i\mid_{(\ref{a1})}=0$ \\
\textbf{Theorem 2.2} \cite{art:2006} Suppose $D_i T^i=0$ is a
conservation law of PDE system (\ref{a1}). Under the contact
transformation, there exist functions $\tilde{T}^i$ such that $J
~D_i T^i=\tilde{D}_i \tilde{T}^i $ where $\tilde{T}^i$ is given
explicitly in terms of the determinant obtained through replacing
the ith row of the Jacobian determinant by$[T^1,T^2,...,T^n]$, where \\
\begin{equation}\label{a21}
J=\left|%
\begin{array}{cccc}
  {\tilde D}_1 x_1&{\tilde D}_1 x_2 & ... &{\tilde D}_1 x_n\\
  {\tilde D}_2 x_1 & {\tilde D}_2 x_2 &...& {\tilde D}_2 x_n \\
  \vdots & \vdots & \vdots &\vdots\\
  {\tilde D}_n x_1 & {\tilde D}_n x_2 &...& {\tilde D}_n x_n \\
\end{array}%
\right|
\end{equation}

\textbf{Theorem 2.3} Suppose $D_i T^i=0$ is a conservation law of
PDE system (\ref{a1}). Under the contact transformation, there
exist functions $\tilde{T}^i$ such that $J ~D_i T^i=\tilde{D}_i
\tilde{T}^i $ where $\tilde{T}^i$ is given explicitly  in terms of
\begin{equation}\label{a22}
\begin{array}{cc}
  \left(%
\begin{array}{cccc}
\tilde{T}^1\\
\tilde{T}^2 \\
\vdots \\
\tilde{T}^n \\
\end{array}%
\right)=J(A^{-1})^T\left(%
\begin{array}{cccc}
  T^1\\
  T^2\\
  \vdots\\
  T^n \\
\end{array}%
\right) , &
J \left(%
\begin{array}{cccc}
   T^1\\
  T^2\\
  \vdots\\
  T^n \\
\end{array}%
\right)= A^T\left(%
\begin{array}{cccc}
  \tilde{T}^1\\
\tilde{T}^2 \\
\vdots \\
\tilde{T}^n \\
\end{array}%
\right) \\
\end{array},
\end{equation}
where
\begin{equation}\label{a23}
\begin{array}{cc}
A=\left(%
\begin{array}{cccc}
  {\tilde D}_1 x_1&{\tilde D}_1 x_2 & ... &{\tilde D}_1 x_n\\
  {\tilde D}_2 x_1 & {\tilde D}_2 x_2 &...& {\tilde D}_2 x_n \\
  \vdots & \vdots & \vdots &\vdots\\
  {\tilde D}_n x_1 & {\tilde D}_n x_2 &...& {\tilde D}_n x_n \\
\end{array}%
\right), &
A^{-1}=\left(%
\begin{array}{cccc}
  {D}_1 \tilde {x_1}&{D}_1 \tilde {x_2} & ... &{ D}_1 \tilde {x_n}\\
  {D}_2 \tilde {x_1}&{D}_2 \tilde {x_2} & ... &{ D}_2 \tilde {x_n}\\
  \vdots & \vdots & \vdots &\vdots\\
  {D}_n \tilde {x_1}&{D}_n \tilde {x_2} & ... &{ D}_n \tilde {x_n}\\
\end{array}%
\right)
\end{array}
\end{equation}
and $J=det(A)$.\\
\textbf{Proof :}\\

Using theorem 2.2 we can write
\begin{equation}
 {\tilde T}^1=\left|%
\begin{array}{cccc}
  T_1&T_2 & ... &T_n \\
  {\tilde D}_2 x_1 & {\tilde D}_2 x_2 &...& {\tilde D}_2 x_n \\
  \vdots & \vdots & \vdots &\vdots\\
  {\tilde D}_n x_1 & {\tilde D}_n x_2 &...& {\tilde D}_n x_n \\
\end{array}%
\right|=
\frac{1}{J}\left|%
\begin{array}{cccc}
  J~T_1& {\tilde D}_2 x_1 & ... &{\tilde D}_n x_1\\
  J~T_2& {\tilde D}_2 x_2 &...&{\tilde D}_n x_2 \\
  \vdots & \vdots & \vdots &\vdots\\
  J~T_n  &{\tilde D}_2 x_n   &...& {\tilde D}_n x_n \\
\end{array}%
\right|,
\end{equation}
\begin{equation}
 {\tilde T}^2=\left|%
\begin{array}{cccc}
   {\tilde D}_1 x_1&{\tilde D}_1 x_2 & ... &{\tilde D}_1 x_n\\
  T_1&T_2 & ... &T_n \\
  \vdots & \vdots & \vdots &\vdots\\
  {\tilde D}_n x_1 & {\tilde D}_n x_2 &...& {\tilde D}_n x_n \\
\end{array}%
\right|=
\frac{1}{J}\left|%
\begin{array}{cccc}
  {\tilde D}_1 x_1 &J~T_1&  ... &{\tilde D}_n x_1\\
  {\tilde D}_1 x_2 &J~T_2& ...&{\tilde D}_n x_2 \\
  \vdots & \vdots & \vdots &\vdots\\
  {\tilde D}_1 x_n   &J~T_n  &...& {\tilde D}_n x_n \\
\end{array}%
\right|,
\end{equation}
\begin{equation}
 {\tilde T}^n=\left|%
\begin{array}{cccc}
   {\tilde D}_1 x_1&{\tilde D}_1 x_2 & ... &{\tilde D}_1 x_n\\
  {\tilde D}_2 x_1 & {\tilde D}_2 x_2 &...& {\tilde D}_2 x_n \\
  \vdots & \vdots & \vdots &\vdots\\
   T_1&T_2 & ... &T_n \\
\end{array}%
\right|=
\frac{1}{J}\left|%
\begin{array}{cccc}
  {\tilde D}_1 x_1& {\tilde D}_2 x_1 & ... &J~T_1\\
  {\tilde D}_1 x_2& {\tilde D}_2 x_2 &...&J~T_2 \\
  \vdots & \vdots & \vdots &\vdots\\
  {\tilde D}_1 x_n &{\tilde D}_2 x_n   &...& J~T_n \\
\end{array}%
\right|.
\end{equation}
Since
\begin{equation}
J=\left|%
\begin{array}{cccc}
  {\tilde D}_1 x_1&{\tilde D}_1 x_2 & ... &{\tilde D}_1 x_n\\
  {\tilde D}_2 x_1 & {\tilde D}_2 x_2 &...& {\tilde D}_2 x_n \\
  \vdots & \vdots & \vdots &\vdots\\
  {\tilde D}_n x_1 & {\tilde D}_n x_2 &...& {\tilde D}_n x_n \\
\end{array}%
\right|=
\left|%
\begin{array}{cccc}
  {\tilde D}_1 x_1& {\tilde D}_2 x_1 & ... &{\tilde D}_n x_1\\
  {\tilde D}_1 x_2& {\tilde D}_2 x_2 &...&{\tilde D}_n x_2 \\
  \vdots & \vdots & \vdots &\vdots\\
  {\tilde D}_1 x_n &{\tilde D}_2 x_n   &...& {\tilde D}_n x_n \\
\end{array}%
\right|=
\left|%
A^T%
\right|,
\end{equation}
one can use the Cramer's rule to find that ${{\tilde T}^1, {\tilde T}^2,\dots,{\tilde T}^n}$ can be written as follows:
\begin{equation}
\left(%
\begin{array}{cccc}
  J~T^1\\
  J~T^2\\
  \vdots\\
  J~T^n \\
\end{array}%
\right)= A^T\left(%
\begin{array}{cccc}
  \tilde{T}^1\\
\tilde{T}^2 \\
\vdots \\
\tilde{T}^n \\
\end{array}%
\right). \\
\end{equation}
Lastly, one can easily see that
\begin{equation}
    A A^{-1}=I.
\end{equation}

\textbf{Lemma 2.1}\\
 Consider $n$ independent variables $x = (x^{1}, x^{2}, ... ,x^{n})$ , $m$ dependent variables $u = (u^{1}, u^{2},..., u^{m})$ and
 the change of independent variables $\tilde{x} = (\tilde{x}^{1}, \tilde{x}^{2}, ... ,\tilde{x}^{n})$, then
  any vector $f(x,u,u_1)=(f^1,f^2,...,f^n)$ must satisfy the
  following identity
\begin{equation}
\left[%
\begin{array}{cccc}
\tilde{D}_1&\tilde{D}_1&\ldots&\tilde{D}_1\\
\tilde{D}_2&\tilde{D}_2 &\ldots&\tilde{D}_2\\
\vdots&\vdots&\ldots&\vdots \\
\tilde{D}_n&\tilde{D}_n &\ldots&\tilde{D}_n\\
\end{array}%
\right]
\left(%
\begin{array}{cccc}
f^1&f^2&\ldots&f^n\\
f^1&f^2&\ldots&f^n\\
\vdots&\vdots&\ldots&\vdots \\
f^1&f^2&\ldots&f^n\\
\end{array}%
\right)=
A \left[%
\begin{array}{cccc}
{D}_1&{D}_1&\ldots&{D}_1\\
{D}_2&{D}_2 &\ldots&{D}_2\\
\vdots&\vdots&\ldots&\vdots \\
{D}_n&{D}_n &\ldots&{D}_n\\
\end{array}%
\right]
\left(%
\begin{array}{cccc}
f^1&f^2&\ldots&f^n\\
f^1&f^2&\ldots&f^n\\
\vdots&\vdots&\ldots&\vdots \\
f^1&f^2&\ldots&f^n\\
\end{array}%
\right),
\end{equation}
where
\begin{equation}
A=\left(%
\begin{array}{cccc}
  {\tilde D}_1 x_1&{\tilde D}_1 x_2 & ... &{\tilde D}_1 x_n\\
  {\tilde D}_2 x_1 & {\tilde D}_2 x_2 &...& {\tilde D}_2 x_n \\
  \vdots & \vdots & \vdots &\vdots\\
  {\tilde D}_n x_1 & {\tilde D}_n x_2 &...& {\tilde D}_n x_n \\
\end{array}%
\right)
\end{equation}
\textbf{Proof :}\\

Since
\begin{equation}
   \tilde{D}_i f^j={\tilde D}_i x_k D_k f^j,~~~ i, j=1,..n, 
\end{equation}
then
\begin{equation}
\left(%
\begin{array}{cccc}
\tilde{D}_1 f^1&\tilde{D}_1 f^2&\ldots&\tilde{D}_1 f^n\\
\tilde{D}_2 f^1&\tilde{D}_2 f^2&\ldots&\tilde{D}_2 f^n\\
\vdots&\vdots&\ldots&\vdots \\
\tilde{D}_n f^1&\tilde{D}_n f^2&\ldots&\tilde{D}_n f^n\\
\end{array}%
\right)= A
\left(%
\begin{array}{cccc}
D_1 f^1&D_1 f^2&\ldots&D_1 f^n\\
D_2 f^1&D_2 f^2&\ldots&D_2 f^n\\
\vdots&\vdots&\ldots&\vdots \\
D_n f^1&D_n f^2&\ldots&D_n f^n\\
\end{array}%
\right)
\end{equation}

{\bf Theorem 2.4} {(Fundamental Theorem of double reduction).}\\
Suppose $D_i T^i=0$ is a conservation law of PDE system
(\ref{a1}). Under the similarity transformation of a symmetry X
of the form (\ref{a5}) for the PDE, there exist functions $\tilde{T}^i$ such
that X is still a symmetry for the PDE $\tilde{D}_i \tilde{T}^i=0
$ and
\begin{equation}
  \left(%
\begin{array}{cccc}
X \tilde{T}^1\\
X \tilde{T}^2 \\
\vdots \\
X \tilde{T}^n \\
\end{array}%
\right)=J(A^{-1})^T\left(%
\begin{array}{cccc}
  $[$T^1,X$]$\\
  $[$T^2,X$]$\\
  \vdots\\
  $[$T^n,X$]$ \\
\end{array}%
\right),
\end{equation}
where
\begin{equation}\label{a23}
\begin{array}{cc}
A=\left(%
\begin{array}{cccc}
  {\tilde D}_1 x_1&{\tilde D}_1 x_2 & ... &{\tilde D}_1 x_n\\
  {\tilde D}_2 x_1 & {\tilde D}_2 x_2 &...& {\tilde D}_2 x_n \\
  \vdots & \vdots & \vdots &\vdots\\
  {\tilde D}_n x_1 & {\tilde D}_n x_2 &...& {\tilde D}_n x_n \\
\end{array}%
\right), &
A^{-1}=\left(%
\begin{array}{cccc}
  {D}_1 \tilde {x_1}&{D}_1 \tilde {x_2} & ... &{ D}_1 \tilde {x_n}\\
  {D}_2 \tilde {x_1}&{D}_2 \tilde {x_2} & ... &{ D}_2 \tilde {x_n}\\
  \vdots & \vdots & \vdots &\vdots\\
  {D}_n \tilde {x_1}&{D}_n \tilde {x_2} & ... &{ D}_n \tilde {x_n}\\
\end{array}%
\right)
\end{array}
\end{equation}
and $J=det(A)$.

\textbf{Proof :}

By the above theorem there exist functions $\tilde{T}^i$ such that
$J ~D_i T^i=\tilde{D}_i \tilde{T}^i $ and
\begin{equation}
\begin{array}{cc}
  \left(%
\begin{array}{cccc}
\tilde{T}^1\\
\tilde{T}^2 \\
\vdots \\
\tilde{T}^n \\
\end{array}%
\right)=J(A^{-1})^T\left(%
\begin{array}{cccc}
  T^1\\
  T^2\\
  \vdots\\
  T^n \\
\end{array}%
\right) , &
J \left(%
\begin{array}{cccc}
   T^1\\
  T^2\\
  \vdots\\
  T^n \\
\end{array}%
\right)= A^T\left(%
\begin{array}{cccc}
  \tilde{T}^1\\
\tilde{T}^2 \\
\vdots \\
\tilde{T}^n \\
\end{array}%
\right) \\
\end{array}
\end{equation}
Then
 X is a symmetry for
the PDE $\tilde{D}_i \tilde{T}^i=0 $, because $X(J) ~D_i T^i+J
~X(D_i T^i)=X(\tilde{D}_i \tilde{T}^i)$ and
\begin{equation}
  \left(%
\begin{array}{cccc}
X \tilde{T}^1\\
X \tilde{T}^2 \\
\vdots \\
X \tilde{T}^n \\
\end{array}%
\right)=J(A^{-1})^T\left(%
\begin{array}{cccc}
  XT^1\\
  XT^2\\
  \vdots\\
  XT^n \\
\end{array}%
\right)+JX((A^{-1})^T)
\left(%
\begin{array}{c}
  T^1 \\
  T^2 \\
 \vdots \\
  T^n \\
\end{array}%
\right)+X(J)(A^{-1})^T
\left(%
\begin{array}{c}
  T^1 \\
  T^2 \\
 \vdots \\
  T^n \\
\end{array}%
\right).
\end{equation}
Since $J=det(A)$, then
\begin{equation}
X(J)= \left|%
\begin{array}{cccc}
  {\tilde D}_1 \xi^1&{\tilde D}_1 \xi^2 & ... &{\tilde D}_1 \xi^n\\
  {\tilde D}_2 x_1 & {\tilde D}_2 x_2 &...& {\tilde D}_2 x_n \\
  \vdots & \vdots & \vdots &\vdots\\
  {\tilde D}_n x_1 & {\tilde D}_n x_2 &...& {\tilde D}_n x_n \\
\end{array}%
\right|
+\left|%
\begin{array}{cccc}
  {\tilde D}_1 x_1&{\tilde D}_1 x_2 & ... &{\tilde D}_1 x_n\\
  {\tilde D}_2 \xi^1 & {\tilde D}_2 \xi^2 &...& {\tilde D}_2 \xi^n \\
  \vdots & \vdots & \vdots &\vdots\\
  {\tilde D}_n x_1 & {\tilde D}_n x_2 &...& {\tilde D}_n x_n \\
\end{array}%
\right|
+ \cdots+\left|%
\begin{array}{cccc}
  {\tilde D}_1 x_1&{\tilde D}_1 x_2 & ... &{\tilde D}_1 x_n\\
  {\tilde D}_2 x_1 & {\tilde D}_2 x_2 &...& {\tilde D}_2 x_n \\
  \vdots & \vdots & \vdots &\vdots\\
  {\tilde D}_n \xi^1 & {\tilde D}_n \xi^2 &...& {\tilde D}_n \xi^n \\
\end{array}%
\right|
\end{equation}
Let $\zeta_{i j}$ denote the cofactor of $ \tilde{D}_i \xi^j$ ,
then it is the cofactor of $ \tilde{D}_i x_j$ for the matrix $A$.
Thus
\begin{equation}
 X(J)=\tilde{D}_i \xi^j ~~\zeta_{i j}=D_k \xi^j ~~\tilde{D}_i x_k~~\zeta_{i j}=D_k \xi^j ~~\delta_{j k}~~J.
\end{equation}
 Since  $\tilde{D}_i x_k~~\zeta_{i j}=\delta_{j k}~~J$ for every fixed $j$ where $\delta_{j k}$ is the Kronecker delta,
then
\begin{equation}
X(J)=J (D_1 \xi^1+D_2 \xi^2 + ... +D_n \xi^n)
\end{equation}
Now using the previous lemma one gets,
\begin{equation}
\left[%
\begin{array}{cccc}
\tilde{D}_1&\tilde{D}_1&\ldots&\tilde{D}_1\\
\tilde{D}_2&\tilde{D}_2 &\ldots&\tilde{D}_2\\
\vdots&\vdots&\ldots&\vdots \\
\tilde{D}_n&\tilde{D}_n &\ldots&\tilde{D}_n\\
\end{array}%
\right]
\left(%
\begin{array}{cccc}
\xi^1&\xi^2&\ldots&\xi^n\\
\xi^1&\xi^2&\ldots&\xi^n\\
\vdots&\vdots&\ldots&\vdots \\
\xi^1&\xi^2&\ldots&\xi^n\\
\end{array}%
\right)=
A \left[%
\begin{array}{cccc}
{D}_1&{D}_1&\ldots&{D}_1\\
{D}_2&{D}_2 &\ldots&{D}_2\\
\vdots&\vdots&\ldots&\vdots \\
{D}_n&{D}_n &\ldots&{D}_n\\
\end{array}%
\right]
\left(%
\begin{array}{cccc}
\xi^1&\xi^2&\ldots&\xi^n\\
\xi^1&\xi^2&\ldots&\xi^n\\
\vdots&\vdots&\ldots&\vdots \\
\xi^1&\xi^2&\ldots&\xi^n\\
\end{array}%
\right)
\end{equation}
Now transposing both sides gives,
\begin{equation}
X(A^T)= \left(%
\begin{array}{cccc}
  D_1 \xi^1&D_2 \xi^1 & ... &D_n \xi^1\\
   D_1 \xi^2&D_2 \xi^2 & ... &D_n \xi^2\\
  \vdots & \vdots & \vdots &\vdots\\
   D_1 \xi^n&D_2 \xi^n & ... &D_n \xi^n\\
\end{array}%
\right)A^T
\end{equation}
 Since $A^T(A^{-1})^T=I$,  then
$X(A^T)(A^{-1})^T+A^T X((A^{-1})^T)=0$, thus
\begin{equation}
\begin{array}{c}
X((A^{-1})^T)=-(A^T)^{-1} X(A^T)(A^{-1})^T=-(A^{-1})^T
X(A^T)(A^T)^{-1}\\\\
=-(A^{-1})^T
\left(%
\begin{array}{cccc}
  D_1 \xi^1&D_2 \xi^1 & ... &D_n \xi^1\\
   D_1 \xi^2&D_2 \xi^2 & ... &D_n \xi^2\\
  \vdots & \vdots & \vdots &\vdots\\
   D_1 \xi^n&D_2 \xi^n & ... &D_n \xi^n\\
\end{array}%
\right)\\
\end{array}.
\end{equation}
Lastly we get the result
\begin{equation}
\left(%
\begin{array}{cccc}
X \tilde{T}^1\\
X \tilde{T}^2 \\
\vdots \\
X \tilde{T}^n \\
\end{array}%
\right)=J(A^{-1})^T
\left(%
\left(%
\begin{array}{cccc}
  XT^1\\
  XT^2\\
  \vdots\\
  XT^n \\
\end{array}%
\right)-\left(%
\begin{array}{cccc}
  D_1 \xi^1&D_2 \xi^1 & ... &D_n \xi^1\\
   D_1 \xi^2&D_2 \xi^2 & ... &D_n \xi^2\\
  \vdots & \vdots & \vdots &\vdots\\
   D_1 \xi^n&D_2 \xi^n & ... &D_n \xi^n\\
\end{array}%
\right)
\left(%
\begin{array}{c}
  T^1 \\
  T^2 \\
 \vdots \\
  T^n \\
\end{array}%
\right) +D_i \xi^i
\left(%
\begin{array}{c}
  T^1 \\
  T^2 \\
 \vdots \\
  T^n \\
\end{array}%
\right)%
\right)
\end{equation}

\textbf{Corollary 2.1} (The necessary and sufficient condition to get reduced conserved form)\\
The conserved form $D_i T^i=0$ of PDE system (\ref{a1}) can be
reduced under the similarity transformation of a symmetry $X$ to a
reduced conserved form $\tilde{D}_i \tilde{T}^i=0$ if and only if
$X$ is associated with the conservation law $T$,
  i.e. $[ T , X ]\mid_{(\ref{a1})}=0 $.

\textbf{Corollary 2.2} (The generalized double reduction theory)\\
A non linear system of $q$ th order PDEs with $n$ independent and
$m$ dependent variables, which admits a nontrivial conserved form
that has at least one associated symmetry in every reduction from
the $n$ reductions (the first step of double reduction) can be
reduced to a non linear system $(q-1)$th order of ODEs .

\textbf{Corollary 2.3} (The inherited symmetries)\\
Any symmetry $Y$ for the conserved form $D_i T^i=0$ of PDE system
(\ref{a1}) can be transformed under the similarity transformation of a
symmetry $X$ for the PDE to the symmetry $\tilde{Y}$ for the PDE
$\tilde{D}_i \tilde{T}^i=0 $.\\

{\bf Remark}:\\
 There is a possibility to get an associated symmetry with a reduced
conserved form by inhering of the non associated symmetry with the
original conserved form. So there is an important useful of the
non associated symmetry also in  Double reduction.\\
Finally we conjecture that the reduction under a combination of an
associated and a non associated symmetries will give us two PDE
one of them is a reduced conserved form and the second is a non
reduced conserved form, we can sperate them via the condition
$X(\tilde{D}_i \tilde{T}^i)=0 $ such that the solution of a
reduced conserved form is also a solution of the non reduced
conserved form.

%%%%%%%%%%%%%%%%%%%%%%%%%%%%%%%%%%%%%%%%%%%%%%%%%%%%%%%%%%%%%%%%%%%%%%%%%%%%%%%%%%%%%%%%%%%%%%%%%%%
\section{Application of the generalized double reduction theory to nonlinear (2 + 1) wave equation }
\setcounter{equation}{0} The  nonlinear (2 + 1) wave equation for
arbitrary function $f(u)$ and $g(u)$
\begin{equation}
    u_{tt}-(f(u)u_x)_x-(g(u)u_y)_y=0,
\end{equation}
has the the obvious conservation law
\begin{equation}
T=\textbf{(}-u_t,f \left( u \right) u_x,g \left( u \right) u_y).\\
\end{equation}
And admits the following four symmetries:

\begin{equation}\label{d1}
\begin{array}{lll}
 X_{1}=\frac{\partial}{\partial t},&X_{2}= \frac{\partial}{\partial x}\\
 X_{3}= \frac{\partial}{\partial y}&X_{4}=t\frac{\partial}{\partial t}+x\frac{\partial}{\partial x}+ y
\frac{\partial}{\partial y}.
\end{array}
\end{equation}
We can get a reduced conserved form for the PDE by the associated
symmetry
 which satisfies the following formula
\begin{equation}
\begin{array}{c}
X\left(%
\begin{array}{c}
  {T}^t\\
  {T}^x\\
  {T}^y\\
\end{array}%
\right)-
\left(%
\begin{array}{ccc}
   D_t \xi^t&D_x \xi^t & D_y \xi^t\\
   D_t \xi^x&D_x \xi^x & D_y \xi^x\\
   D_t \xi^y&D_x \xi^y & D_y \xi^y\\
\end{array}%
\right)
\left(%
\begin{array}{c}
  {T}^t\\
  {T}^x\\
  {T}^y\\
\end{array}%
\right)\\
 +(D_t \xi^t+D_x \xi^x+D_y \xi^y)
\left(%
\begin{array}{c}
  {T}^t\\
  {T}^x\\
  {T}^y\\
\end{array}%
\right)=0.
\end{array}
\end{equation}
Then the only associated symmetries are $X_{1},X_{2}$ and $X_{3}$,
so we can get a reduced conserved form by the combination of them
$X=\frac{\partial}{\partial t}+c_1\frac{\partial}{\partial
x}+c_2\frac{\partial}{\partial y}$, where the generator $X$ has a
canonical form $X=\frac{\partial}{\partial q}$ when
\begin{equation}
    \frac{dt}{1}= \frac{dx}{c_1}= \frac{dy}{c_2}=\frac{du}{0}=\frac{dr}{0}=\frac{ds}{0}=\frac{dq}{1}=
    \frac{dw}{0},
\end{equation}
or
\begin{equation}
   \begin{array}{cccc}
          r = y-c_2 t, & s = x-c_1 t, & q =t,& w(r, s) = u. \\
    \end{array}
\end{equation}
Using the following formula, we can get the reduced
conserved form
\begin{equation}
\left(%
\begin{array}{c}
  {T}^r\\
  {T}^s\\
  {T}^q\\
\end{array}%
\right)=J(A^{-1})^T\left(%
\begin{array}{c}
  {T}^t\\
  {T}^x\\
  {T}^y\\
\end{array}%
\right),
\end{equation}
where
\begin{equation}
\begin{array}{cc}
A^{-1}=\left(%
\begin{array}{ccc}
  {D}_t r&{D}_t s & { D}_t q \\
  {D}_x r&{D}_x s & { D}_x q \\
  {D}_y r&{D}_y s & { D}_y q \\
\end{array}%
\right), &J=det(A).
\end{array}
\end{equation}
Then the reduced conserved form is
\begin{equation}
 D_r T^{r}+ D_s T^{s}=0,
\end{equation}
where
\begin{equation}
\begin{array}{lll}
 T^{r}={c_2}^2 w_r+c_2 c_1 w_s-g(w) w_r,\\
 T^{s}=c_1 c_2 w_r+c_1^2w_s-f(w)w_s, \\
 T^{q} =-c_2w_r-c_1w_s.
\end{array}
\end{equation}
The reduced conserved form admits the inherited symmetry: \be
\begin{array}{l}
 \tilde{X}_{4}=r \frac{\partial}{\partial r}+s \frac{\partial}{\partial s},\\
\end{array}\ee

Similarly we can get a reduced conserved form for the PDE by the
associated symmetry
 which satisfies the following formula
\begin{equation}
\begin{array}{c}
X\left(%
\begin{array}{c}
  {T}^r\\
  {T}^s\\
\end{array}%
\right)-
\left(%
\begin{array}{ccc}
   D_r \xi^r&D_s \xi^r \\
   D_r \xi^s&D_s \xi^s \\
\end{array}%
\right)
\left(%
\begin{array}{c}
  {T}^r\\
  {T}^s\\
\end{array}%
\right)\\
 +(D_r \xi^r+D_s \xi^s)
\left(%
\begin{array}{c}
  {T}^r\\
  {T}^s\\
\end{array}%
\right)=0.
\end{array}
\end{equation}
One can see that $\tilde{X}_{4}$ is an associated symmetry,
so we can get a reduced conserved form by  $Y=r\frac{\partial}{\partial r}+ s\frac{\partial}{\partial s}$, where
the generator $Y$ has a canonical form $Y=\frac{\partial}{\partial
m}$ when
\begin{equation}
    \frac{dr}{ r}= \frac{ds}{ s}= \frac{dw}{0}=\frac{dn}{0}=\frac{dm}{1}=
    \frac{dv}{0},
\end{equation}
or
\begin{equation}
   \begin{array}{ccc}
          n=\frac{s}{r}, & m=ln r ,& v(n)=w.\\
       \end{array}
\end{equation}
So by using the following formula, we can get the reduced
conserved form
\begin{equation}
\left(%
\begin{array}{ccccc}
  {T}^n\\
  {T}^m\\
\end{array}%
\right)=J(A^{-1})^T\left(%
\begin{array}{c}
 {T}^r \\
 {T}^s\\
\end{array}%
\right),
\end{equation}
where
\begin{equation}
\begin{array}{cc}
A^{-1}=\left(%
\begin{array}{ccc}
  {D}_r n&{D}_r m  \\
  {D}_s n&{D}_s m \\
\end{array}%
\right), &J=det(A).
\end{array}
\end{equation}
Then the reduced conserved form is:
\begin{equation}
 D_n T^{n}=0,
\end{equation}
where
\begin{equation}
\begin{array}{lll}
 T^{n}= v_n(-{c_2}^2 n^2+2c_2 c_1 n+n^2g(v)-{c_1}^2+f(v)),\\
 T^{m}=-v_n(-c_2^2 n+c_2c_1+n g(v)). \\
\end{array}
\end{equation}
The second step of double reduction can be given as
\begin{equation}
  v_n(-{c_2}^2 n^2+2c_2 c_1 n+n^2g(v)-{c_1}^2+f(v))=C,
\end{equation}
where $C$ is a constant, $n=\frac{x-c_1 t}{y-c_2 t} $ and  $v=u.$
%%%%%%%%%%%%%%%%%%%%%%%%%%%%%%%%%%%%%%%%%%%%%%%%%%%%%%%%%%%%%%%%%%%%%%%%%%%%%%%%%%%%%%%%%%%%%%%%%%%
\section{Conclusion}
We have shown that the double reduction theory is still true in general case. This shows that one
can obtain \emph{the invariant solution for a non linear system of
PDEs by this procedure from the association of the symmetry with
its conserved form via the new generalized formula (\ref {a22}).}
%%%%%%%%%%%%%%%%%%%%%%%%%%%%%%%%%%%%%%%%%%%%%%%%%%%%%%%%%%%%%%%%%%%%%%%%%%%%%%%%%%%%%%%%%%%%%%%%%%%%%%%%%%%%%%%%%%%%%%%%%%%%%%%%%%%%%%%%%%%%%%%%%%%%%%%%%%%%%%%%%%%%%%%%%%%%%%%%%%%%%%%%%%%%%%%%%%%%%%

\end{document}